\documentclass[a4paper,10pt]{amsart}
\usepackage{amsfonts}
\usepackage{amsthm}
\usepackage{amssymb}
\usepackage{amsmath}
\usepackage{enumerate}
\usepackage{url}
\usepackage{color}

\begin{document}
\title[On large values of $L(\sigma,\chi)$]{On large values of $L(\sigma,\chi)$}

\author[C. Aistleitner]{Christoph Aistleitner}
\address{Institute of Analysis and Number Theory, TU Graz, Steyrergasse 30, 8010 Graz, Austria}
\email{aistleitner@math.tugraz.at}

\author[K. Mahatab]{Kamalakshya Mahatab}
\address{Department of  Mathematical Sciences, Norwegian University of Science and Technology, NO-7491 Trondheim, Norway}
\email{accessing.infinity@gmail.com}

\author[M. Munsch]{Marc Munsch}
\address{Institute of Analysis and Number Theory, TU Graz, Steyrergasse 30, 8010 Graz, Austria}
\email{munsch@math.tugraz.at}

\author[A. Peyrot]{Alexandre Peyrot}
\address{Department of Mathematics, Stanford University, 450 Serra Mall, Building 380, Stanford CA 94305-2125}
\email{alpeyrot@stanford.edu}

\newcommand{\mods}[1]{\,(\mathrm{mod}\,{#1})}

\begin{abstract}
In recent years a variant of the resonance method was developed which allowed to obtain improved $\Omega$-results for the Riemann zeta function along vertical lines in the critical strip. In the present paper we show how this method can be adapted to prove the existence of large values of $|L(\sigma, \chi)|$ in the range $\sigma \in (1/2,1]$, and to estimate the proportion of characters for which $|L(\sigma, \chi)|$ is of such a large order. More precisely, for every fixed $\sigma \in (1/2,1)$ we show that for all sufficiently large $q$ there is a non-principal character $\chi \mods q$ such that $\log |L(\sigma,\chi)| \geq C(\sigma) (\log q)^{1-\sigma} (\log \log q)^{-\sigma}$. In the case $\sigma=1$ we show that there is a non-principal character $\chi \mods q$ for which $|L(1,\chi)| \geq e^\gamma \left(\log_2 q + \log_3 q - C \right)$. In both cases, our results essentially match the prediction for the actual order of such extreme values, based on probabilistic models. 
\end{abstract}

\maketitle

\def \l {{\lambda}}
\def \a {{\alpha}}
\def \b {{\beta}}
\def \f {{\phi}}
\def \r {{\rho}}
\def \R {{\mathbb R}}
\def \H {{\mathbb H}}
\def \N {{\mathbb N}}
\def \C {{\mathbb C}}
\def \Z {{\mathbb Z}}
\def \F {{\Phi}}
\def \Q {{\mathbb Q}}
\def \e {{\epsilon }}
\def\GL{\ensuremath{\mathop{\textrm{\normalfont GL}}}}
\def\SL{\ensuremath{\mathop{\textrm{\normalfont SL}}}}
\def\Gal{\ensuremath{\mathop{\textrm{\normalfont Gal}}}}
\def\SU{\ensuremath{\mathop{\textrm{\normalfont SU}}}}
\def\SO{\ensuremath{\mathop{\textrm{\normalfont SO}}}}

\newtheorem{prop}{Proposition}
\newtheorem{claim}{Claim}
\newtheorem{lemma}{Lemma}
\newtheorem{thm}{Theorem}
\newtheorem{defn}{Definition}
\newtheorem{conj}{Conjecture}

\theoremstyle{definition}
\newtheorem{exmp}{Example}

\theoremstyle{remark}
\newtheorem{rmk}{Remark}

\section{Introduction and statement of results}

A classical problem in analytic number theory is to understand the values of Dirichlet $L$-series in the critical strip. The study of these values and their distribution is a deep question in number theory which has various important repercussions for the related attached arithmetic, algebraic and geometric objects. In the last century, the notion of family of $L$-functions has been important  both as a heuristic guide to understand or guess many important statistical properties of $L$-functions.  Since the initial discussion between Montgomery and Dyson about correlations of zeros of the Riemann zeta function, a large amount of work has since been carried out in order to model families of $L$-functions by characteristic polynomials of random matrices. As addressed by Katz and Sarnak \cite{katz} in a general way, we can naturally attach a group of matrices to every family of $L$-functions which assigns a ``symmetry type" to it. For instance, in a concrete treatment, the conjectural estimation of moments of $L$-functions in various families of $L$-functions regarding  their type of symmetry is comprehensively discussed in \cite{con}. Two important families of $L$-functions are known to have the same symmetry coming from an unitary group of matrices: \begin{itemize}
\item $\left\{L(\sigma +it),  \,t \geq 0\right\}$, ordered by $t$, where $L(s)$ is any primitive $L$-function.
\item $\left\{L(s, \chi), \,\chi \text{ a primitive character modulo }  q\right\}$ ordered by $q$. \end{itemize}

Thus it can be expected, at least on a heuristic level, to get similar statistical results for these two families. The most famous and oldest example of an $L$-function is the Riemann zeta function, which is naturally related to the distribution of prime numbers. The importance of large values of the Riemann zeta function and $L$-functions is known since Littlewood, and studies of this sort of problem originated in the close connections with historical number theoretical problems about character sums, quadratic non-residues as well as class numbers.\\

There are different methods to prove the existence of large values of the Riemann zeta function along vertical lines in the critical strip. Montgomery \cite{mont} used a method based on Diophantine approximation to prove that
\begin{equation}\label{Montgomery} 
\max_{t \in [T,2T]} \log\vert \zeta(\sigma+it)\vert \geq C(\sigma) (\log T)^{1-\sigma} (\log_2 T)^{-\sigma}.
\end{equation}
Here and throughout this paper we write $\log_j$ for the $j$-th iterated logarithm, so for example $\log_2 q = \log \log q$. Balasubramanian and Ramachandra \cite{bala} obtained comparable results, including unconditional results for the case $\sigma=1/2$, using large moments of the zeta function and the connection with divisor sums. A more recent method is the so-called resonance method, whose main functional principle was probably first used by Voronin \cite{voronin}. However, the main development of the method was due to Soundararajan \cite{sound}, who realized how it leads to a general maximization problem for quadratic forms, and who showed the potential strength as well as the wide potential range of applications of the method. Hilberdink \cite{hilb} took up Soundararajan's ideas and made important contributions to the refinement of the method. Basically, the idea of this method is to find a function $R(t)$ such that the quotient $I_1/I_2$ is ``large'', where
$$
I_1 = \int_0^T \zeta (\sigma + it) |R(t)|^2 ~dt
$$
and
$$
I_2=\int_0^T |R(t)|^2 ~dt.
$$
A lower bound for the quotient $I_1/I_2$ immediately yields a corresponding lower bound for the maximal value of $|\zeta(\sigma+it)|$ in the range $t \in [0,T]$. This method, of simple appearance, turned out to be quite powerful. Recently, using ``long'' resonators as they were introduced in \cite{aist}, Bondarenko and Seip \cite{bs} achieved a breakthrough by obtained improved $\Omega$-results for the Riemann zeta function on the critical line. In a recent work involving three of the the authors \cite{amm}, a ``long'' resonator of completely multiplicative type was used to prove that there exists a constant $C$ such that 
\begin{equation} \label{1-line}
 \max_{t\in [\sqrt{T},T]}|\zeta(1+it)|\geq e^{\gamma}(\log_2 T +\log_3 T -C).
\end{equation}
We note that the in the cases of $\sigma \in (1/2,1)$ and $\sigma=1$ the known lower bounds essentially match the conjectured truth (up to the values of the constants involved), while for $\sigma=1/2$ the situation remains unclear even on a conjectural level (cf.\ also \cite{farm}).\\

The main purpose of the present paper is to show that it is possible to adapt the resonance method in such a way that it allows to prove the existence of non-principal characters $\chi \mods q$ for which $|L(\sigma,\chi)|$ is large. For this purpose the quotient $I_1/I_2$ from above is replaced by the quotient of finite sums
$$
\frac{\left|\sum_{\chi \mods q} L(\sigma,\chi) |R(\chi)^2| \right|}{\sum_{\chi \mods q} |R(\chi)|^2},
$$
where $R(\chi)$ is a function which ``resonates'' with $L(\sigma,\chi)$. Our proofs will show that it is possible to define the function $R(\chi)$ as a short Euler product in a completely multiplicative way. As for the Riemann zeta function, the two cases $\sigma=1$ and $\sigma \in (1/2,1)$ differ significantly in several respects, and we will first discuss one and then the other case. We do not address the case $\sigma=1/2$ in the present paper, since in that case it is better to use a resonator which is not completely multiplicative, which requires a different approach from the one taken in our paper (cf.\ also the remark at the end of the introduction).\\ 

In the case $\sigma=1$, it is known that
\begin{equation} \label{lowup}
q^{-\varepsilon} \ll_\varepsilon |L(1,\chi)| \ll \log q,
\end{equation}
and that assuming the Generalized Riemann Hypothesis this can be improved to
\begin{equation} \label{little}
\frac{1}{\log_2 q} \ll |L(1,\chi)| \ll \log_2 q.
\end{equation}
Much effort has been spent on possible improvements of the upper bound in \eqref{lowup}, where, however, anything beyond improvements of the implied constant seems to be completely out of reach. See for example \cite{gs1,lou,ram}. On the other hand Little\-wood's result \eqref{little} is essentially optimal, up to the value of the implied constant, since it is known that there are characters for which $|L(1,\chi)|$ is of order roughly $1/(\log_2 q)$ and $\log_2 q$, respectively. The strongest results concerning the existence of ``extremal'' characters are due to Granville and Soundararajan \cite{gs2}, who proved that for sufficiently large $q$ and any given $A \geq 10$ there are at least $q^{1-1/A}$ characters $\chi$ (mod $q$) for which 
\begin{equation} \label{gslarge}
|L(1,\chi)| \geq e^\gamma (\log_2 q + \log_3 q - \log_4 q - \log A - C),
\end{equation}
for some absolute constant $C$. The same paper also contains results for small values of $|L(1,\chi)|$, which establish the existence of characters $\chi$ for which $1/|L(1,\chi)| \geq \frac{6 e^\gamma}{\pi^2} (\log_2 q - C)$. In Theorem \ref{th1} below we improve Granville and Soundararajan's result to a level which matches the corresponding result for the Riemann zeta function noted in \eqref{1-line}, and which also essentially matches the prediction for the optimal result for this problem (up to the values of the constants involved). The subsequent Theorem \ref{th2} contains an estimate for the proportion of characters for which $|L(1,\chi)|$ is of this order of magnitude, and also is in general accordance with the prediction for the actual truth based on probabilistic models.\\

\begin{thm} \label{th1}
Let $\varepsilon >0$. For all sufficiently large $q$ there is a non-principal character $\chi \mods q$ such that
$$
|L(1,\chi)| \geq e^\gamma \left(\log_2 q + \log_3 q - C - \varepsilon \right),
$$
where 
$$
C = 1 + \log_2 4 \approx 1.33.
$$
\end{thm}

\begin{thm} \label{th2}
Let $C$ denote the constant from Theorem \ref{th1}. Let $\delta>0$, and define 
\begin{equation} \label{quant}
\Phi(\delta) = \Big|\big\{\chi \mods q : ~|L(1,\chi)| > e^\gamma (\log_2 q + \log_3 q - C - \delta) \big\}\Big|.
\end{equation}
Then, as $q$ tends to infinity, we have
$$
\Phi(\delta) \geq q^{1-e^{-\delta}+o(1)}.
$$
\end{thm}

Theorem \ref{th2} should be compared with estimates for the proportion of characters with large values of $|L(1,\chi)|$ due to Granville and Soundararajan \cite{gs2}. They proved that the number of characters $\chi \mods q$ for which $|L(1,\chi)| > e^\gamma \tau$ is 
$$
q \cdot \exp\left( -\frac{2 e^{\tau-C_0-1}}{\tau} (1 + \text{errors}) \right)
$$
for some explicit constant $C_0 \approx -0.395$, valid up to $\tau \leq \log_2 q - 20$. If this asymptotic formula could be extended to values of $\tau$ corresponding to extremely large values of $|L(1,\chi)|$ such as those in Theorem \ref{th1}, then for all sufficiently large $q$ one could deduce the existence of a non-principal character $\chi \mods q$ for which
\begin{equation} \label{conj1}
|L(1,\chi)| \geq e^\gamma \left(\log_2 q + \log_3 q + C_0 + 1 - \log 2 + o(1)\right),
\end{equation}
together with the density estimate
\begin{align}
& \Big|\big\{\chi \mods q : ~|L(1,\chi)| > e^\gamma (\log_2 q + \log_3 q + C_0 + 1 - \log 2 - \delta) \big\}\Big| \nonumber\\ 
& = q^{1 - e^{-\delta} + o(1)}.  \label{conj2}
\end{align}
Based on a comparison with the distribution of values of random Euler products, Granville and Soundararajan speculated that these formulas might show the actual truth about extremal values of $|L(1,\chi)|$ for characters $\chi \mods q$. We remark that the constant $-C \approx -1.33$ in our theorems does not match with the constant $C_0 + 1 - \log 2 \approx -0.09$ in these hypothetical estimates, and that our density estimate has $e^{-\delta}$ while \eqref{conj2} would give the better factor $e^{-\delta-(C_0+1-\log 2 + C)} \approx 0.29 e^{-\delta}$ to estimate the quantity in \eqref{quant}. Still, generally speaking our results are in good accordance with the conjecture of Granville and Soundararajan, and if \eqref{conj1} and \eqref{conj2} actually represent the truth, then our results are essentially optimal up to the precise values of the constants.\\

We note that the method which is used for the proof of Theorem \ref{th2} can be adapted in such a way that it gives a lower bound for the measure of those $t \in [0,T]$ for which $|\zeta(1+it)|$ is of size around $e^\gamma (\log_2 T + \log_3 T)$. More precisely, it can be shown that
$$
\textup{meas} \Big(t \in [0,T]:~ |\zeta(1+it)| \geq e^\gamma (\log_2 T + \log_3 T - C - \delta) \Big) \geq T^{1 - e^{-\delta} + o(1)},
$$
where $C$ is the same constant as in Theorems \ref{th1} and \ref{th2}. This improves results of Granville and Soundararajan \cite{gs2}, where an additional term $- \log_4 T$ is necessary.\\

Now we turn to the case $\sigma \in (1/2,1)$. A paper of Lamzouri \cite{lam2} provides very precise estimates for the distribution of large values of various families of $L$-functions at a fixed point $1/2<\sigma<1$, strenghtening the range of validity obtained in \cite{gs3}. As noted in \cite{lam2}, as a by-product Lamzouri obtained results on the existence of extreme values inside a family of $L$-functions, ordered by conductor. More precisely, he could deduce that in the corresponding family there exists an $L$-function of order $\exp\left( C(\sigma) (\log Q)^{1-\sigma}/\log_2 Q\right)$, where $Q$ is the conductor of the family. These results are of the same quality as the ones that could be obtained as a consequence of Soundararajan resonance method \cite{sound}. Though, as pointed out by Lamzouri, this method fails to provide estimates inside the critical strip as strong as the $\Omega$-results obtained by Montgomery for the Riemann zeta function. On the other hand, a drawback of Montgomery's method is that it does not generalize to other families of $L$-functions. Thus it is certainly interesting to prove such a result for any other family of $L$-functions. Assuming the Generalized Riemann Hypothesis (GRH), Lamzouri was able to achieve this for the family of $L$-functions associated to Legendre symbols averaged over primes $p\leq x$. More precisely, in \cite[Theorem $6$]{lam2} he proved that there exists a constant $C(\sigma)>0$ such that for $\delta>0$  sufficiently small and $x>0$ large enough, 
\begin{equation}\label{Legendre} \log |L(s,\chi_p)|\geq (C(\sigma)+o(1))(\log x)^{1-\sigma}(\log_2 x)^{-\sigma}\end{equation} holds for more than $\delta x^{1/2}$ primes $p\leq x$. Theorem \ref{th3} below gives a Montgomery-type $\Omega$-result in the family attached to Dirichlet characters $\mods q$, for $q$ sufficiently large. As in Theorem \ref{th1}, our result essentially matches the conjecturally optimal result, except for the precise value of the constant.\\

\begin{thm}\label{th3}  Let $1/2<\sigma<1$. There exists a constant $C(\sigma)>0$ such that for all sufficiently large $q$ there exists a non-principal character $\chi \mods q$ such that
\begin{equation}\label{large} 
\log \vert L(\sigma,\chi)\vert \geq C(\sigma)(\log q)^{1-\sigma}(\log_2 q)^{-\sigma}.
\end{equation}
\end{thm}

We note that Lamzouri averages over a specific subset of well chosen primes under GRH. Basically this comes down to fixing the first values of the character $\left(\frac{.}{p}\right)$, an idea which traces back to Littlewood \cite{little} and Chowla \cite{chowla}, and was subsequently used by Montgomery \cite{mont2} and Granville and Soundararajan \cite{gs3}. Lamzouri used GRH again to approximate $L(\sigma,\chi)$ in this relatively ``thin" set of characters and to control precisely the off-diagonal terms arising from the orthogonality of quadratic characters. Since we work unconditionally and use zero-density estimates, we cannot restrict ourselves to such a thin family. \\

We close this section with some remarks on remaining open problems, and possible improvements of our theorems. As in other implementations of the ``long resonator'', a certain positivity property plays a crucial role in the argument. This explains why our method cannot be applied for example to the case of  $1/|L(\sigma,\chi)|$, since the positivity property is lost in this setting. Secondly, it would be interesting to adapt our method to prove the existence of large values of $|L(1,\chi_d)|$ for primitive real characters $\chi_d$, as $d$ varies over all fundamental discriminants in a range $|d| \leq x$; cf.\ \cite{gs3,lam}. However, this seems to be very difficult to implement, since the orthogonality relations between these characters are much more subtle than those in our setting, and again one looses the required positivity property.\\

It is worth pointing that the constant $C(\sigma)$ appearing in Theorem \ref{th3} comes both from the particular choice of the resonator as well as from the zero-density estimates used in the proof. It might be possible to optimize these constants by analyzing carefully the zero-density estimates and changing slightly the resonator chosen. We decided to not carry out this task in order to keep the presentation short. Furthermore, the method which is used in Theorem \ref{th2} to estimate the proportion of characters for which $|L(1,\chi)|$ is large could be applied in the setting of Theorem \ref{th3} as well; however, again to keep the paper short, but also since the results which one would obtain in this case are significantly weaker than those in Theorem \ref{th2}, we have restrained from carrying this out. In the case $\sigma \in (1/2,1)$, a key requirement of our method is to control the contribution of ``bad" characters. A similar argument would easily give the existence of large values for even characters, or any specified character lying in a subgroup $H$ of sufficiently large order compared to the order of the full group of Dirichlet characters.\\

\emph{Remark:} \quad While this paper was under revision, a new paper of de la Bret\`eche and Tenenbaum \cite{dlb} appeared on the arxiv. Their paper contains a comprehensive analysis of the resonance method on the critical line as well as several new applications. Among their results are improvements of the $\Omega$-results of Bondarenko and Seip for the Riemann zeta function, a complete solution of the related problem concerning the maximal order of GCD sums, applications to large values of character sums, as well as results for extreme values of $|L(1/2,\chi)|$ which complement those obtained in the present paper. Furthermore, there have been other new applications of the resonance method to related problems. These applications cover the argument of the zeta function \cite{bs_arg,chir}, Hardy's $Z$-function \cite{kamal}, and short character sums \cite{munsch}.

\section{Proof of Theorem \ref{th1}} \label{sec:th1}

Let $q$ be ``large'', choose $B>\log 4$, and define $X = \log q \log_2 q / B$. We set  
$$
L(1,\chi,X) = \prod_{p \leq X} \left(1 - \frac{\chi(p)}{p} \right)^{-1} =: \sum_{k=1}^\infty a_k \chi(k),
$$
where $a_k$ can only be either $0$ or $1/k$. For primes $p \leq X$ we set
$$
q_p := \left(1-\frac{p}{X}\right),
$$
while for primes $p > X$ we set $q_p = 0$, and we extend this to $q_m$ for general $m$ in a completely multiplicative way. We write $G_q$ for the group of Dirichlet characters $\mods q$. For a given character $\chi \in G_q$, we define the resonator function $R(\chi)$ as
$$
R(\chi) := \prod_{p\leq X} (1-q_p\chi(p))^{-1}.
$$


Finally, let $Y = \exp \left((\log q)^{20} \right)$ and set 
$$
L(1,\chi,Y) = \prod_{p \leq Y} \left(1 - \frac{\chi(p)}{p} \right)^{-1} =: \sum_{k=1}^\infty b_k \chi(k),
$$
where again $b_k$ can only be either $0$ or $1/k$. Note that $b_k \geq a_k$ for all $k \geq 1$, since $Y \geq X$.\\

With this choice of $Y$ it is not difficult to establish the approximation formula
\begin{equation} \label{lapprox}
L(1,\chi) = L(1,\chi,Y) \left(1 + \mathcal{O} \left(\frac{1}{(\log q)^2} \right) \right),
\end{equation}
which holds for all characters $\chi \mods q$ except for $\chi_0$ and at most one other exceptional character, which we denote by $\chi_{e}$ (note that is widely believed that such an exceptional character does not actually exist). To establish this approximation formula, we can argue as in the proof of the prime number theorem (similar to the sketch which is given at the beginning of the proof of Theorem 2 in \cite{gs2} for a finite approximation to the Riemann zeta function). More precisely, we note that by Perron's formula
$$
\int_{2+i \infty}^{2-i \infty} \log L(1+s,\chi) \frac{Y^s}{s} ds = -\sum_{p \leq Y} \log \left( 1 - \frac{\chi(p)}{p} \right) + \underbrace{\mathcal{O} \left( \sum_{p \leq Y} \sum_{k \geq \frac{\log Y}{\log p}} \frac{1}{kp^k}\right)}_{= \mathcal{O}\left(\frac{1}{\log Y}\right)}.
$$
Now we truncate the integral at $T = Y^3$ and shift the contour to the left of 0, in such a way that we enclose a zero-free region of $L(1,\chi)$. By a classical result going back to Gronwall, Landau and Titchmarsh there is at most one exceptional character $\chi_{e} \mods q$ such that for all other characters $\chi \neq \chi_{e}$ we have $L(\sigma+it,\chi) \neq 0$ for $\sigma \geq 1 - \frac{c}{\log (q(2+|t|))}$ for some constant $c$ (see for example \cite{hb}). Thus for $\chi \neq \chi_0,\chi_{e}$ the only pole of the function in the integral which is enclosed by the curve is at $s=0$, and it is a standard procedure to establish upper bounds for all error terms and to derive the desired approximation formula.\\

Set 
$$
S_1 = \sum_{\chi \in G_q} L(1,\chi,Y) |R(\chi)|^2 
$$
and 
$$
S_2 = \sum_{\chi \in G_q} |R(\chi)|^2.
$$
By expanding $|R(\chi)|^2$ and the two finite approximations $L(1,\chi,X) = \sum_k a_k \chi(k)$ and $L(1,\chi,Y)=\sum_{k} b_k \chi(k)$, and by using the orthogonality of characters together with the facts that $b_k \geq a_k$ for all $k$ and that $(a_n)_n,(b_n)_n,(q_n)_n$ are all non-negative, we have
\begin{align}
S_1 & = \sum_{\chi \in G_q} \sum_k \sum_{m,n} b_k \chi(k) q_m q_n \chi(m) \overline{\chi(n)}  \nonumber \\
& = \sum_k ~ b_k \phi(q) \sum_{\substack{m,n:\\~km \equiv n \mods q}} q_m q_n \nonumber \\
& \geq \sum_k ~ a_k \phi(q)\sum_{\substack{m,n:\\~km \equiv n \mods q}} q_m q_n. \label{akqm}
\end{align}
The advantage of changing from $L(1,\chi,Y)$ to $L(1,\chi,X)$ is that the coefficients $(a_n)_n$ are supported on the same set as $(q_n)_n$, namely on the set of $X$-smooth numbers. This will allow us to obtain a fully explicit formula for the quotient $S_1/S_2$ in equation \eqref{expl} below.\\

By expanding $S_2$ we have
\begin{align}
S_2 & = \sum_{\chi \in G_q} \sum_{m,n} q_mq_n\chi(m)\overline{\chi(n)} \nonumber\\
& = \phi(q) \sum_{\substack{m,n:\\ m \equiv n \mods q}} q_mq_n. \label{i2rep}
\end{align}

Now assume that $k$ is fixed. Exploiting the fact that the resonator $R(\chi)$ is constructed in a completely multiplicative way we have 
\begin{align}
a_k \phi(q) \sum_{\substack{m,n:\\~km \equiv n \mods q}}  q_m q_n & \geq a_k \phi(q) \sum_{\substack{m,n: ~k|n,\\~km \equiv n \mods q}} q_m q_n \nonumber\\ 
& = a_k \phi(q) \sum_{\substack{r,m:\\ km \equiv kr \mods q}} q_m \underbrace{q_{kr}}_{=q_k q_r} \nonumber\\
& \geq a_k q_k ~\underbrace{\phi(q) \sum_{\substack{r,m:\\r \equiv m \mods q}} q_r q_m}_{= S_2 \text{~by \eqref{i2rep}}}. \label{relation}
\end{align}
Note that again the fact that $(q_n)_n$ are non-negative numbers played a crucial role in this calculation.\\

Thus by \eqref{akqm} and \eqref{relation} we have 
\begin{align}
\frac{S_1}{S_2} & \geq \sum_{k=1}^\infty a_k q_k \nonumber \\
& = \prod_{p \leq X} \left(1 - q_p p^{-1} \right)^{-1} \label{expl} \\
& = \left( \prod_{p \leq X} (1 - p^{-1})^{-1} \right) \left( \prod_{p \leq X} \frac{p-1}{p-q_p} \right). \label{two_prod}
\end{align}
For the first product in line \eqref{two_prod} we have
\begin{equation} \label{mert}
\prod_{p \leq X} (1 - p^{-1})^{-1} \geq e^\gamma \log X \left(1 - \frac{1}{2(\log X)^2}\right)
\end{equation}
by Mertens' third theorem with an explicit estimate for the error term from \cite[Theorem 8]{rs}. For the other product in \eqref{two_prod} we have
\begin{align}
- \log \left( \prod_{p \leq X} \frac{p-1}{p-q_p} \right) & = - \sum_{p \leq X} \log \left(1 - \frac{p}{p+(p-1)X}\right) \nonumber \\
& \leq\sum_{p \leq X} \left(\frac{1}{X} + \frac{2}{pX} \right) \nonumber\\
& \leq \left(1+\mathcal{O}\left(\frac{1}{(\log X)^2} \right)\right) \frac{1}{\log X} \label{otherpb}
\end{align}
by the prime number theorem. Thus overall we have
\begin{equation} \label{i1i2ratio}
\frac{S_1}{S_2} \geq e^\gamma \log X \left(1 - \frac{1}{\log X} +\mathcal{O}\left(\frac{1}{(\log X)^2} \right) \right).
\end{equation}

We note that 
\begin{align} 
\log |R(\chi_0)|^2 & = \log\prod_{p\leq X} (1-q_p)^{-2} \nonumber\\
& = 2 \sum_{p \leq X} (\log X - \log p) \nonumber\\
& = 2(1+o(1))\frac{X}{\log X},\label{boundr}
\end{align}
and that by the prime number theorem and by our choice of $Y$
\begin{equation} \label{boundl}
L(1,\chi_0,Y) \ll (\log q)^{20}.
\end{equation}

To obtain a lower bound for $S_2$ we note that by \eqref{i2rep}
\begin{align}
\nonumber S_2 & \geq \phi(q) \sum_m q_m^2\\
& = \phi(q) \prod_{p\leq X} (1-q_p^2)^{-1}.
\end{align}
We now estimate by partial summation
\begin{align*}
\log\prod_{p\leq X} (1-q_p^2)^{-1} & = - \sum_{p\leq X} \log\left(\frac{2pX-p^2}{X^2}\right)\\
&=2  \int_2^X \frac{\pi(x) (X-x)}{2xX-x^2} \mathrm{d}x\\
&= 2(1+o(1)) \int_2^X \frac{X-x}{(\log x)(2X-x)}\mathrm{d}x\\
&= 2(1+o(1)) \left(\frac{X}{\log X} - X\int_2^X \frac{1}{(\log x)(2X-x)}\mathrm{d}x\right)\\
&=  2(1+o(1)) \left(\frac{X}{\log X} - X \int_X^{2X-2} \frac{1}{(\log(2X-u))u} \mathrm{d}u\right)\\
& = 2(1+o(1)) \left(\frac{X}{\log X} - X \int_1^{2-2/X} \frac{1}{(\log(2-t) + \log X) t} \mathrm{d}t\right).
\end{align*}
We now let 
$$J:= \int_1^{2-2/X} \frac{1}{(\log(2-t) + \log X) t} \mathrm{d}t\ = J_1 + J_2,$$
where $J_1$ corresponds to the contribution of the integral from $1$ to $2-(\log X)^{-2}$ and $J_2$ to the contribution of the integral from  $2-(\log X)^{-2}$ to $2-2/X$. One easily checks that
$$J_1 = (1+o(1)) \frac{\log 2}{\log X},$$
while $J_2 = \mathcal{O}((\log X)^{-2})$, so that 
$$\log\prod_{p\leq X} (1-q_p^2)^{-1}  = (2-\log 4) (1+o(1)) \frac{X}{\log X}.$$
Thus in total we have
\begin{align} 
S_2 & \geq \phi(q) \exp \left((1+o(1)) (2-\log 4) \frac{X}{\log X}\right) \nonumber\\
& = \exp \left( (1+o(1)) \left( \left(1+\frac{2-\log 4}{B} \right) \log q \right) \right),\label{order1}
\end{align}
and the same lower bound holds for $S_1$ due to \eqref{i1i2ratio}.\\

On the other hand, by \eqref{boundr} and \eqref{boundl}
\begin{equation} \label{order2}
|R(\chi_0)|^2 \leq \exp\big( (1+o(1)) (2\log q)/B \big)
\end{equation}
and similarly
\begin{equation} \label{order3}
| L(1,\chi_0,Y) R(\chi_0)^2| \leq \exp\big( (1+o(1)) (2 \log q)/B \big),
\end{equation}
and the same upper bounds hold for $R(\chi_{e})$ instead of $R(\chi_0)$. Note that we assumed $B > \log 4$, which guarantees that $1+(2- \log 4)/B > 2/B$, such that the expressions in \eqref{order2} and \eqref{order3} are asymptotically much smaller than that in \eqref{order1}. This allows us to remove the contribution of $\chi_0$ and $\chi_{e}$ from the sums in $S_1$ and $S_2$, and together with \eqref{mert} and \eqref{otherpb} we obtain
\begin{align}
\frac{\left|{\sum}^* ~L(1,\chi,Y) |R(\chi)|^2 \right|}{{\sum}^* ~|R(\chi)|^2} & \geq e^\gamma \log X \left(1 - \frac{1}{\log X} + \mathcal{O}\left(\frac{1}{(\log X)^2} \right)\right) \nonumber\\
& = e^\gamma \left(\log_2 q + \log_3 q - \log B - 1 + \mathcal{O}\left(\frac{1}{\log_2 q} \right)\right), \label{i1i2obt}
\end{align}
where the asterisk indicates that both summations extend over all characters $\chi \in G_q$ under the restriction that $\chi \neq \chi_0,\chi_{e}$. 
Together with \eqref{lapprox} this implies that for sufficiently large $q$ there is a non-principal character for which
$$
|L(1,\chi)| \geq e^\gamma \left(\log_2 q + \log_3 q - \log B - 1 + o\left(1\right) \right),
$$
which proves the theorem since we can choose $B$ as close to $\log 4$ as we wish.

\section{Proof of Theorem \ref{th2}}

Let $\delta > 0$ be given. Set $B = e^\delta e^{-\frac{1}{\sqrt{\log_2 q}}} \log 4$. Then $B$ is bounded away from $\log 4$ (for sufficiently large $q$), so we can use the construction and notations from the proof of Theorem \ref{th1}. Set $C = 1 + \log_2 4$ and 
\begin{align*}
A_\delta & = e^\gamma \left(\log_2 q + \log_3 q - C - \delta + \frac{1}{2 \sqrt{\log_2 q}} \right),\\
A_{\delta}' & = e^\gamma \left(\log_2 q + \log_3 q - C - \delta \right).
\end{align*}
Note that the additional factor $e^{-\frac{1}{\sqrt{\log_2 q}}}$ in the definition of $B$ leads to an extra $+\frac{1}{\sqrt{\log_2 q}}$ in \eqref{i1i2obt} which dominates asymptotically over the error term $\mathcal{O}\left(\frac{1}{\log_2 q} \right)$, so that for sufficiently large $q$ the right-hand side of \eqref{i1i2obt} exceeds $A_\delta$. Furthermore, the approximation formula \eqref{lapprox} allows us to switch from $L(1,\chi,Y)$ to $L(1,\chi)$ in \eqref{i1i2obt}, and again the error can be absorbed in the difference between $A_\delta$ and the right-hand side of \eqref{i1i2obt}. Thus, writing again ${\sum}^*$ for the sum over all characters $\mods q$ except for $\chi_0$ and the potential exceptional character $\chi_{e}$, from \eqref{i1i2obt} we deduce that for sufficiently large $q$ we have
\begin{align*}
A_\delta {\sum}^*~ |R(\chi)|^2 & \leq {\sum}^*~ L(1,\chi) |R(\chi)|^2\\\\
&\leq \left|~\sideset{}{^*}\sum_{\chi: ~|L(1,\chi)| \leq A_\delta'} L(1,\chi)|R(\chi)|^2 \right| +  \left|~ \sideset{}{^*}\sum_{\chi: ~|L(1,\chi)| > A_\delta'} L(1,\chi)|R(\chi)|^2\right|\\
&\leq A_\delta' {\sum}^* |R(\chi)|^2  +  \left|~ \sideset{}{^*}\sum_{\chi: ~|L(1,\chi)| > A_\delta'} L(1,\chi)|R(\chi)|^2 \right|.
\end{align*}
We thus obtain 
\begin{equation}\label{keyeq}
\left|~ \sideset{}{^*}\sum_{\chi: ~|L(1,\chi)| > A_\delta'} L(1,\chi) |R(\chi)|^2 \right| \geq \frac{e^\gamma}{2 \sqrt{\log_2 q}} {\sum}^* |R(\chi)|^2.
\end{equation}
By a result of Granville and Soundararajan \cite{gs1} we have 
\begin{equation}\label{maxGS}
\max_{\chi \not=\chi_0}|L(1,\chi)| \leq (1+o(1)) \frac{\log q}{3}.
\end{equation}
Furthermore, we have
$$
\left|~ \sideset{}{^*}\sum_{\chi: ~|L(1,\chi)| > A_\delta'} L(1,\chi) |R(\chi)|^2 \right|\leq \max_{\chi \not= \chi_0} |R(\chi)|^2 \max_{\chi \not= \chi_0}|L(1,\chi)| \Phi(\delta),
$$
where $\Phi(\delta)$ is the number of characters mod $q$ for which $|L(1,\chi)|$ exceeds $A_{\delta}'$, as defined in the statement of Theorem \ref{th2}. Consequently, by (\ref{keyeq}) and (\ref{maxGS}),
\begin{align*}
\Phi(\delta) & \geq (1+o(1)) \frac{3 e^\gamma {\sum}^* |R(\chi)|^2}{2(\log q)(\log_2 q)^{1/2} \max |R(\chi)|^2}.
\end{align*}
By \eqref{boundr} we have
\begin{equation*}\label{maxR}
\max|R(\chi)|^2 \leq |R(\chi_0)|^2 \leq e^{(2+o(1))X/\log X} \leq q^{(2+o(1))/B},
\end{equation*}
while by \eqref{order1} and the subsequent remarks we have
$$
{\sum}^* |R(\chi)|^2 \geq q^{(1+o(1)) \left(1 + \frac{2-\log 4}{B} \right)}.
$$
Putting all these estimates together, we finally obtain
$$
\Phi(\delta) \geq q^{(1+o(1)) \left(1- \frac{\log 4}{B} \right)},
$$
which by our choice of $B$ implies
$$
\Phi(\delta) \geq q^{1-e^{-\delta}+o(1)}.
$$
This proves Theorem \ref{th2}.

\section{Proof of Theorem \ref{th3}}

The proof of Theorem \ref{th3} follows roughly the same path as the proof of Theorem \ref{th1}. However, there are two major differences. On the one hand, the conclusion of Theorem \ref{th1} is very precise and it is important to avoid large errors. For Theorem \ref{th2} we do not have to be so careful, and for example we can simply ignore the contribution of the prime powers since it is negligible in comparison with the main term. On the other hand, the approximation of the $L$-function by a finite Euler product is more complicated in the case $\sigma \in (1/2,1)$, since only weaker results on the zeros of the $L$-functions can be used. Thus instead of considering two exceptional characters as in the proof of Theorem \ref{th1}, now we have to take care of a potentially much larger class of ``bad'' exceptional characters.\\

We recall the following approximation lemma for Dirichlet $L$-functions. This type of result, when combined with powerful zero-density estimates, shows that with very few exceptions the $L$-function, when averaged over a suitable family, can be approximated by very short Euler products (over the primes $p\leq (\log Q)^A$, where $Q$ is the conductor of the corresponding family) in the strip $1/2<\textup{Re}(s) <1$. 

\begin{lemma}\cite[Lemma 8.2]{JAMS}\label{approxDir} Let $q$ be  a large prime and let $\chi$ be a character $(\text{mod } q)$. Let  $X \geq 2$ and $|t|\leq 3q$ be real numbers.
Let $\frac{1}{2} \leq \sigma_0 <\sigma\leq 1$ and suppose that the
rectangle $\{ s: \ \ \sigma_0 <\textup{Re}(s) \leq 1, \ \
|\textup{Im}(s) -t| \leq X+2\}$ does not contain any zeros of $L(s,\chi)$.
Then
$$
\log L(\sigma+it,\chi)= \sum_{n=2}^{X} \frac{\Lambda(n)\chi(n)}{n^{\sigma+it}\log n} +
\mathcal{O}\left( \frac{\log
q}{(\sigma_1-\sigma_0)^2}X^{\sigma_1-\sigma}\right),
$$
where  $\sigma_1 = \min(\sigma_0+\frac{1}{\log X},
\frac{\sigma+\sigma_0}{2})$.
\end{lemma}

If we assume the Generalized Riemann Hypothesis for $L(s,\chi)$, then the hypotheses of Lemma \ref{approxDir} are verified with the parameter $\sigma_0=1/2$, which allows to approximate our $L$-function by a short Euler product. However, even without assuming GRH, the lemma remains very useful since the zero-density estimates imply that we have good finite approximations for most characters $\chi$. More precisely, let $N(\sigma,H, \chi)$ denote the number of zeros of $L(s,\chi)$ such that $\textup{Re}(s)\geq \sigma$ and $|\textup{Im}(s)|\leq H$. It follows directly from a zero-density result of Montgomery \cite[Theorem 12.1]{mont2} (Equations ($12.9$) and ($12.10$)) that for $H \geq 2$ and $1/2<\sigma<1$  we have
$ \sum_{\chi \in G_q} N(\sigma,H, \chi) \ll (qH)^{3(1-\sigma)/(2-\sigma)}(\log qH)^{14}$. Thus, combining this estimate and Lemma \ref{approxDir} with $t=0$,
 $X=(\log q)^{3/(\sigma-1/2)}$, $H=X+2$ and $\sigma_0=\sigma/2+1/4>1/2$, we obtain
\begin{equation}\label{approxsum} \log L(\sigma,\chi)=\sum_{n=2}^X\frac{\Lambda(n)\chi(n)}{n^{\sigma}\log n}+\mathcal{O}\left(\frac{1}{(\log q)^{1/4}}\right),\end{equation}
for all characters $\chi \ (\text{mod }q)$ except for a set of  ``bad" characters $\text{Bad}_{\sigma}(q)$ of cardinality $\leq q^{1-a(\sigma)}$ for some constant $a(\sigma)>0$.\\

As in Section \ref{sec:th1}, we let $G_q$ denote the group of Dirichlet characters modulo $q$. We assume that $\sigma$ in the range $1/2<\sigma<1$ is fixed, and set $X=(\log q)^{3/(\sigma-1/2)}$. For every $\chi \in G_q$ we define a Dirichlet polynomial $D_\chi(\sigma,X)$ associated to $L(\sigma,\chi)$ by setting
$$
D_\chi(\sigma,X) = \sum_{2 \leq n \leq X} \frac{\Lambda(n)\chi(n)}{n^{\sigma}\log n}.
$$ 
As shown in Lemma \ref{approxDir} above, for our particular choice of $X=(\log q)^{3/(\sigma-1/2)}$ this truncation is a good approximation of $ \log L(\sigma,\chi)$, provided that the character $\chi$ is not contained in the exceptional set $\text{Bad}_{\sigma}(q)$. As we are only interested in non-trivial characters, we set $\text{Bad}^*(q)= \text{Bad}(q) \cup \{\chi_0\}$, where $\chi_0$ is the trivial character modulo $q$. Thus we have 

$$\log L(\sigma,\chi) = D_\chi(\sigma,X) + \mathcal{O}\left(\frac{1}{(\log q)^{1/4}}\right), \chi \notin \text{Bad}^*_{\sigma}(q).$$  Since the bound we are aiming for is not affected by multiplicative error terms, we restrict ourselves to a sum over primes. Indeed, the terms when $n=p^k$ for $k\geq 2$ only contribute $\ll 1$ in total to $D_{\chi}(\sigma,X)$. It follows that
\begin{equation} \label{S:approx}
\log L(\sigma,\chi) = S_{\chi}(\sigma,X) + \mathcal{O}(1), \hspace{2mm}\chi \notin \text{Bad}^*_{\sigma}(q),
\end{equation}
where $\displaystyle{S_{\chi}(\sigma,X)= \sum_{p \leq X} \frac{\chi(p)}{p^{\sigma}}}.$ Consequently, in order to prove Theorem \ref{th3}, it suffices to exhibit large values of $S_{\chi}(\sigma,X)$ for $\chi \notin \text{Bad}^*_{\sigma}(q)$. To do so, we apply a ``long'' resonantor to the approximation $S_{\chi}(\sigma,X)$ in a similar way as performed in Section \ref{sec:th1}. Since we are not aiming for a result which is as precise as Theorem \ref{th1}, we can use a simpler resonator, which still is of completely multiplicative type. For $a(\sigma)$ as in the previous section, we set $Y=\frac{a(\sigma)}{2}(\log q)(\log\log q)$. Furthermore, we set $q_1=1$ and $q_p=0$ for $p>Y$, and $q_p=1/2$ for all small primes $p\leq Y$. We extend this in a completely multiplicative way to obtain weights $q_n$ for all $n\geq 1$. We now define for $\chi\in G_q$ 
\begin{equation} \label{res:def}
R(\chi)=\prod_{p\leq Y} \left(1-q_p\chi(p)\right)^{-1} = \prod_{p\leq Y} \left(1-\frac{\chi(p)}{2} \right)^{-1},
\end{equation}
and note that we can write $R(\chi)$ as a Dirichlet series in the form
\begin{equation} \label{rdi}
R(\chi) = \sum_{n=1}^\infty q_n \chi(n).
\end{equation}
Accordingly we have
$$
|R(\chi)|^2 = \sum_{m,n = 1}^\infty q_m q_n \chi(m)\overline{\chi(n)}.
$$ 
Similar to the previous sections, we consider the sums
$$
S_1 = \sum_{\chi \in G_q} S_{\chi}(\sigma,X)\vert R(\chi)\vert^2 
$$ 
and 
$$
S_2=\sum_{\chi \in G_q} \vert R(\chi)\vert^2,
$$ 
and we will use the simple inequality
\begin{equation}\label{resonance} 
\max_{\chi \in G_q} \vert S_{\chi}(\sigma,X)\vert \geq \frac{ \vert S_1 \vert}{S_2 }. 
\end{equation} 
Note that additionally to a lower bound for the right-hand side of \eqref{resonance}, we also need to show that the contribution to $S_1$ and $S_2$ coming from the characters in $\text{Bad}^*_{\sigma}(q)$ is harmless. Since the set $\text{Bad}^*_{\sigma}(q)$ is relatively small, it is enough to bound individually $R(\chi)$ for $\chi \in G_q$. From the definition \eqref{res:def} and our choice of $Y$ we immediately obtain
\begin{equation}\label{R_est}
|R(\chi)|^2 \leq 2^{2 \pi(Y)} \leq \exp(2Y/\log Y) \leq q^{a(\sigma) \log 2 + o(1)}.
\end{equation}

Summing over bad characters, we get 
\begin{equation}\label{S2bad}
\sum_{\chi \in \text{Bad}_{\sigma}^*(q)}  |R(\chi)|^2  \ll q^{1-a(\sigma)}q^{a(\sigma) \log 2+o(1)} \ll q^{1-\frac{a(\sigma)}{4}+o(1)}.
\end{equation} 
In the other hand, we have trivially that \begin{equation}\label{S2tout} \sum_{\chi \in G_q} \vert R(\chi)\vert^2 =\sum_{\chi\in G_q} q_m q_n \chi(m)\overline{\chi(n)}=\phi(q) \left(\sum_{m=n \bmod q} q_m q_n\right) \geq q^{1-o(1)}\end{equation} where we used together the fact that $q_1=1$ and $\phi(q) \geq q/\log \log q$. Similarly, we can disregard the contribution of the elements of $\text{Bad}_{\sigma}^*(q)$ to $S_1$. Indeed,

$$|S_{\chi}(\sigma,X)| \leq \sum_{p=2}^{X}\frac{1}{p^{\sigma}} \ll X^{1-\sigma} = (\log q)^{\frac{6(1-\sigma)}{2\sigma-1}}=q^{o(1)}.$$ Thus, as before in (\ref{S2bad}) we get
\begin{equation}\label{S1bad}
\left| \sum_{\chi \in \text{Bad}_{\sigma}^*(q)} S_{\chi}(\sigma,X) |R(\chi)|^2 \right|  \ll q^{1-\frac{a(\sigma)}{4}+o(1)}.
\end{equation} 
Expanding $|R(\chi)|^2 $ and switching the order of summation, we have 
\begin{eqnarray*}
S_1& = & \sum_{\chi \in G_q}\vert R(\chi)\vert^2 S_{\chi}(\sigma,X)\\ 
& = & \sum_{p=1}^{X} \frac{1}{p^{\sigma}} \left( \sum_{m,n = 1}^\infty q_m q_n 
 \sum_{\chi \in G_q}\chi(m)\overline{\chi(n)}\chi(p)\right),
\end{eqnarray*} where the inner sum is positive from the orthogonality relations on $G_q$. Thus,

\begin{equation}\label{quotientdev} S_{1} = \sum_{p=1}^{X}\frac{1}{p^{\sigma}} \left(\phi(q)\sum_{\substack{m,n, pm=n \bmod q}} q_m q_n \right).\end{equation} Assume $p$ to be fixed such that $(p,q)=1$. Then, using the positivity and the completely multiplicative property of the coefficients $q_n$ we get 

\begin{eqnarray*}
\phi(q)\sum_{\substack{m,n,\\ pm=n \bmod q}}  q_m q_n  &\geq& \phi(q)\ \sum_{\substack{m,n,  p | n, \\pm=n \bmod q}}  q_m q_n  \\ & = &\phi(q) \sum_{\substack{m,l, pm=pl \bmod q}} q_m \underbrace{q_{pl}}_{=q_p q_l}  = q_p \phi(q) \left(\sum_{m=l \bmod q} q_m q_l\right) \\
&=& q_p S_2. \\
\end{eqnarray*} Noticing that $q_p = 1/2$ for $p \leq Y$ and $q_p=0$ for $p>Y$, we deduce from (\ref{quotientdev}) that

\begin{eqnarray}\label{quotientfinal}\frac{S_1}{S_2} &\geq & \sum_{\substack{p=2,\\ (p,q)=1}}^{X} \frac{q_p}{p^{\sigma}} = \sum_{p=2}^{Y}\frac{1}{2p^{\sigma}} \nonumber \\
& \gg_\sigma & \frac{Y^{1-\sigma}}{\log Y}.\end{eqnarray} Thus, we get

$$\max_{\chi \in G_q} |S_{\chi}(\sigma,X)| \gg_{\sigma} (\log q)^{1-\sigma} (\log \log q)^{-\sigma}.$$ Appealing to (\ref{S2bad}), (\ref{S2tout}) and (\ref{S1bad}), the inequality (\ref{quotientfinal}) remains true when we sum only over characters $\chi \notin \text{Bad}_{\sigma}^*(q)$ yielding to

$$ \max_{\chi \notin \text{Bad}_{\sigma}^*(q)} |S_{\chi}(\sigma,X)| \gg_{\sigma} (\log q)^{1-\sigma} (\log \log q)^{-\sigma}.$$ On this set, $S_{\chi}(\sigma,X)$ is a good approximation to $\log L(\sigma,\chi)$, as noted in \eqref{S:approx}. This proves Theorem \ref{th3}.

\section*{Acknowledgements}

CA is supported by the Austrian Science Fund (FWF), projects Y-901 and F 5512-N26. KM is supported by Grant 227768 of the Research Council of Norway. A part of this work was carried out while he held an ERCIM ``Alain Bensoussan'' Fellowship. MM is supported by FWF project Y-901. AP is supported by the Swiss National Science Foundation grant P2ELP2\_172089. He was supported by the Simons Investigators Grant of Kannan Soundararajan during his time at Stanford University, when this paper was written.


\end{document}